\documentclass[12pt,a4paper]{amsart}
\usepackage{amssymb,color,currfile,pdfpages}

\usepackage[english]{babel}
\usepackage{verbatim,here}
\usepackage[T1]{fontenc}
\usepackage{epstopdf}
\usepackage{floatflt,graphicx}
\usepackage{color,xcolor}
\usepackage{enumerate}
\usepackage{animate}
\usepackage{a4wide}
\usepackage{amsmath, amssymb}
\usepackage{amsthm}
\usepackage{float}
\usepackage{pdfpages}
\usepackage{orcidlink}
\usepackage{hyperref}

\newcounter{minutes}
\divide\time by 60
\newcounter{hours}
\multiply\time by 60 \addtocounter{minutes}{-\time}

\newenvironment{pf}[1][]{%
 \vskip 3mm
 \noindent
 \ifthenelse{\equal{#1}{}}%
  {{\slshape Proof. }}%
  {{\slshape #1.} }%
 }%
{\qed\bigskip}

\dedicatory{}
\commby{}
\theoremstyle{plain}
\newtheorem{thm}[equation]{Theorem}
\newtheorem{cor}[equation]{Corollary}
\newtheorem{lem}[equation]{Lemma}
\newtheorem{example}[equation]{Example}
\newtheorem{prop}[equation]{Proposition}

\theoremstyle{definition}

\theoremstyle{remark}
\newtheorem{rem}[equation]{Remark}
\newtheorem{conj}[equation]{Conjecture}

\numberwithin{equation}{section}

\newcommand{\beq}{\begin{equation}}
\newcommand{\eeq}{\end{equation}}
\newcommand{\ben}{\begin{enumerate}}
\newcommand{\een}{\end{enumerate}}
\newcommand{\bequu}{\begin{eqnarray*}}
\newcommand{\eequu}{\end{eqnarray*}}
\newcommand{\bequ}{\begin{eqnarray}}
\newcommand{\eequ}{\end{eqnarray}}
\newcommand{\B}{\mathbb{B}^2}

\newcommand{\Bn}{ {\mathbb{B}^n} }

\begin{document}
\thispagestyle{empty}
\def\thefootnote{}
\newcommand{\re}{{\rm Re\,}}
\newcommand{\capacity}{{\mathop{\mathrm{cap}}}}
\newcommand{\capacitydenL}{{\mathop{\mathrm{cap\ \underline{dens}}}}}
\newcommand{\capacitydenU}{{\mathop{\mathrm{cap\ \overline{dens}}}}}

\newcommand{\co}{{\overline{\operatorname{co}}}}
\newcommand{\T}{{\mathcal T}}
\newcommand{\U}{{\mathcal U}}
\newcommand{\es}{{\mathcal S}}
\newcommand{\LU}{{\mathcal{LU}}}
\newcommand{\ZF}{{\mathcal{ZF}}}
\newcommand{\IR}{{\mathbb R}}
\newcommand{\IN}{{\mathbb N}}
\newcommand{\IC}{{\mathbb C}}
\newcommand{\IT}{{\mathbb T}}
\newcommand{\ID}{{\mathbb D}}
\newcommand{\IB}{{\mathbb B}}
\newcommand{\K}{{\mathcal K}}
\newcommand{\X}{{\mathcal X}}
\newcommand{\PP}{{\mathcal P}}
\newcommand{\uhp}{{\mathbb H}}
\newcommand{\Z}{{\mathbb Z}}
\newcommand{\N}{{\mathcal N}}
\newcommand{\M}{{\mathcal M}}
\newcommand{\SCC}{{\mathcal{SCC}}}
\newcommand{\CC}{{\mathcal C}}
\newcommand{\st}{{\mathcal{SS}}}
\newcommand{\D}{{\mathbb D}}
\newcommand{\sphere}{{\widehat{\mathbb C}}}
\newcommand{\image}{{\operatorname{Im}\,}}
\newcommand{\Aut}{{\operatorname{Aut}}}
\newcommand{\real}{{\operatorname{Re}\,}}
\newcommand{\kernel}{{\operatorname{Ker}}}
\newcommand{\ord}{{\operatorname{ord}}}
\newcommand{\id}{{\operatorname{id}}}
\newcommand{\mob}{{\text{\rm M\"{o}b}}}
\newcommand{\Int}{{\operatorname{Int}\,}}
\newcommand{\Sign}{{\operatorname{Sign}}}
\newcommand{\inv}{^{-1}}
\newcommand{\area}{{\operatorname{Area}}}
\newcommand{\eit}{{e^{i\theta}}}
\newcommand{\ucv}{{\operatorname{UCV}}}
%%%%%%%%%%%%%%%%%%%%%%%%%%%%%%%%%%%%%%%%%%%%%%%%%%%%%%%%%%%%%%%%%%%%%%%%%%%%%%%%%5

%\newcommand{\pad}[2]{\frac{\der #1}{\der #2}}
\def\be{\begin{equation}}
\def\ee{\end{equation}}
\newcommand{\sep}{\itemsep -0.01in}
\newcommand{\seps}{\itemsep -0.02in}
\newcommand{\sepss}{\itemsep -0.03in}
\newcommand{\bee}{\begin{enumerate}}
\newcommand{\eee}{\end{enumerate}}
\newcommand{\pays}{\!\!\!\!}
\newcommand{\pay}{\!\!\!}
\newcommand{\blem}{\begin{lem}}
\newcommand{\elem}{\end{lem}}
\newcommand{\bthm}{\begin{thm}}
\newcommand{\ethm}{\end{thm}}
\newcommand{\bcor}{\begin{cor}}
\newcommand{\ecor}{\end{cor}}
\newcommand{\beg}{\begin{example}}
\newcommand{\eeg}{\end{example}}
\newcommand{\begs}{\begin{examples}}
\newcommand{\eegs}{\end{examples}}
\newcommand{\bdefe}{\begin{defin}}
\newcommand{\edefe}{\end{defin}}
\newcommand{\bprob}{\begin{prob}}
\newcommand{\eprob}{\end{prob}}
\newcommand{\bei}{\begin{itemize}}
\newcommand{\eei}{\end{itemize}}

\newcommand{\bcon}{\begin{conj}}
\newcommand{\econ}{\end{conj}}
\newcommand{\bcons}{\begin{conjs}}
\newcommand{\econs}{\end{conjs}}
\newcommand{\bprop}{\begin{propo}}
\newcommand{\eprop}{\end{propo}}
\newcommand{\br}{\begin{rem}}
\newcommand{\er}{\end{rem}}
\newcommand{\brs}{\begin{rems}}
\newcommand{\ers}{\end{rems}}
\newcommand{\bo}{\begin{obser}}
\newcommand{\eo}{\end{obser}}
\newcommand{\bos}{\begin{obsers}}
\newcommand{\eos}{\end{obsers}}
\newcommand{\bpf}{\begin{pf}}
\newcommand{\epf}{\end{pf}}
\newcommand{\ba}{\begin{array}}
\newcommand{\ea}{\end{array}}
\newcommand{\llra}{\longleftrightarrow}
\newcommand{\lra}{\longrightarrow}
\newcommand{\lla}{\longleftarrow}
\newcommand{\Llra}{\Longleftrightarrow}
\newcommand{\Lra}{\Longrightarrow}
\newcommand{\Lla}{\Longleftarrow}
\newcommand{\Ra}{\Rightarrow}
\newcommand{\La}{\Leftarrow}
\newcommand{\ra}{\rightarrow}
\newcommand{\la}{\leftarrow}
\newcommand{\ds}{\displaystyle}
\newcommand{\psubset}{\subsetneq}

\def\cc{\setcounter{equation}{0}   % THIS CLEARS THE COUNTER
\setcounter{figure}{0}\setcounter{table}{0}}

\def\cc{\setcounter{equation}{0}   % THIS CLEARS THE COUNTER
\setcounter{figure}{0}\setcounter{table}{0}}
%\DeclareMathOperator*{\esssup}{ess\,sup}
%\DeclareMathOperator*{\essinf}{ess\,inf}
%\DeclareMathOperator*{\Har}{\mathrm{Har}}

%=====================================================================
%\newcounter{minutes}\setcounter{minutes}{\time}
%\divide\time by 60
%\newcounter{hours}\setcounter{hours}{\time}
%\multiply\time by 60 \addtocounter{minutes}{-\time}
%=====================================================================

\title[Angular limits]
{On angular limits of quasiregular mappings}

\author[J. Huang]{Jie Huang \orcidlink{0009-0000-4566-6947}}
\address{ Department of Mathematics with Computer Science, Guangdong Technion, 241 Daxue Road, Jinping District, Shantou, Guangdong 515063, People's Republic of China and Department of Mathematics, Technion - Israel Institute of Technology, Haifa 32000, Israel\newline
\href{https://orcid.org/0009-0000-4566-6947}{{\tt https://orcid.org/0009-0000-4566-6947}}
}
\email{jie.huang@gtiit.edu.cn}

\author[A. Rasila]{Antti Rasila$^{*}$  \orcidlink{0000-0003-3797-942X}}
\address{Department of Mathematics and Statistics,
Department of Mathematics with Computer Science, Guangdong Technion, 241 Daxue Road, Jinping District, Shantou, Guangdong 515063, People's Republic of China and Department of Mathematics, Technion - Israel Institute of Technology, Haifa 32000, Israel\newline
\href{https://orcid.org/0000-0003-3797-942X}{{\tt https://orcid.org/0000-0003-3797-942X}}
}
\email{antti.rasila@iki.fi; antti.rasila@gtiit.edu.cn}

\author[M. Vuorinen]{Matti Vuorinen \orcidlink{0000-0002-1734-8228}}
\address{Department of Mathematics and Statistics,
University of Turku, Finland\newline
\href{https://orcid.org/0000-0002-1734-8228}{{\tt https://orcid.org/0000-0002-1734-8228}}
}
\email{vuorinen@utu.fi}

\date{\today}

%%%%%%%%%%%%%%%%%%%%%%%%%%%%%%%%%%%%%
\begin{abstract}
We investigate Lindel\"of and Koebe type boundary behavior results for bounded quasiregular mappings in $n$-dimensional Euclidean space.
Our results give sufficient conditions for the existence of non-tangential limits at a boundary point.
\end{abstract}

\keywords{Quasiregular mappings, boundary behavior, capacity density, hyperbolic metrics. \\${}^{\mathbf{*}}$ Corresponding author}

\maketitle

\footnotetext{\texttt{{\tiny File:~\jobname .tex, printed: \number\year-%
\number\month-\number\day, \thehours.\ifnum\theminutes<10{0}\fi\theminutes}}}

\section{Introduction}

Recall the following classical result of Lindel\"of (see e.g. 
\cite[p. 259]{Rudin:1987}):

\begin{thm}
\label{linde}
Suppose that $\gamma$ is a curve, with parameter
interval $[0,1]$, such that $|\gamma(t)|<1$ if $t<1$ and $\gamma(1)=1$.
If $f$ is a bounded analytic function of the unit disk $\D$ and
$$
\lim_{t\to 1} f\big(\gamma(t)\big) = \alpha,
$$
then $f$ has angular limit $\alpha$ at $1$, i.e. limit in each angular
region contained in the unit disk with the vertex in $1$.
\end{thm}

It is a natural question to ask if generalizations of Lindel\"of's theorem 
can be obtained for other classes of functions. In particular, as quasiregular 
and quasiconformal maps in $\IR^n$
are respectively natural generalizations of analytic and conformal functions 
of one complex variable, it is interesting to study problems of this type for 
these classes of functions. 

Indeed, it is known that Lindel\"of's theorem holds for quasiconformal 
mappings in any dimension $n\geq 2$. The result in $\IR^3$ was proved 
by Gehring \cite[p. 21]{Gehring}, and the same proof applies to the 
general $n$-dimensional setting as well. Improvements of this result 
have been obtained e.g. in \cite{Vu1}, see also \cite[Chapter 15]{VuorinenBook}.

Interestingly, Lindel\"of's theorem does not hold for non-univalent 
bounded quasiregular mappings for $n\geq 3$, which was shown by the 
counterexample of Rickman (see \cite{Ri1}
or \cite{Ri2}). Some results in the opposite direction
for bounded quasiregular 
mappings $f: \mathbb{B}^n \to  \IR^n$ of the unit ball $\mathbb{B}^n,n \ge 2,$ having a limit $\alpha$ at 
a boundary point $b \in \partial  \mathbb{B}^n$ along a
non-tangential set $E$ were obtained in \cite{Vu2}. Indeed, if the lower capacity
density of $E$ at $b$ is positive, then $f$ has an angular limit $\alpha$ at $b$, see Theorem \ref{Theorem 3.1} below. 
The proof was based on a two constants theorem for quasiregular mappings 
\cite[p. 190]{Ri2}
(for latest publications on the two constants theorem see p. 15 of \cite{HKV}).
The same method was applied to prove further results on the existence 
of angular limits  by the second  and the third authors 
and their collaborators. See e.g. \cite{VuorinenBook,YangRasila}, 
and references therein.

In this paper, we continue this investigation by proving three new results. 
The first of our main results, Theorem \ref{main1}, weakens the assumption of 
\cite[Theorem 3.1]{Vu2} by replacing the assumption that the mapping attains 
a limit along a set of positive lower capacity density at the boundary point 
by a weaker assumption and under the additional requirement that the function
approaches the limit value fast enough. Theorem \ref{main2} further extends the result 
giving a stronger conclusion guaranteeing that the mapping actually has a 
limit along a certain tangential set at the boundary point. Finally, in 
Section \ref{koebe}, we will investigate so-called Koebe type boundary 
results. The aim here is to study results that give sufficient conditions 
for the mapping to be constant, in terms of its behavior on the set that 
has an accumulation point at the boundary.

\section{Preliminary results}

In this section, we recall basic definitions and preliminaries from the literature.

\subsection*{Notations} We shall follow standard notations and terminology adopted from \cite{HKV,vaisala,VuorinenBook}. For $x\in \IR^n$, $n\geq 2$, and for $r>0$ let $B^n(x, r):=\{z\in\IR^n: |z-x|<r\}$, and $S^{n-1}(x, r)=\partial B^n(x, r)$. For the convenience of writing, we set $B^n(r)=B^n(0, r)$, $S^{n-1}(r)=\partial B^n(r)$, $\Bn=B^n(1)$ and $S^{n-1}=\partial \Bn$. 
The surface area of $S^{n-1}$ is denoted by $\omega_{n-1}$ and $\Omega_n$ is the volume of $\Bn$. It is well known that $\omega_{n-1}=n\Omega_n$ and that
$$
\Omega_n=\frac{\pi^{n/2}}{\Gamma(1+n/2)}
$$
for $n=2,3,\ldots$, where $\Gamma$ is Euler's gamma function.

For $b\in S^{n-1}$ and $\varphi\in(0, \pi/2)$, let $C(b,\varphi)$ be the open bounded cone
$$C(b,\varphi):=\big\{x: \langle b,  b-x\rangle >|b-x|\cos(\varphi)\big\}\cap B^n\big( b, \cos(\varphi)\big),$$
where $ \langle x,y  \rangle $  is the inner product $ \langle x,y  \rangle  :=\sum\limits_{i=1}^n x_iy_i$.
If $\cos(\varphi)>r>s>0$, then we write
\begin{equation}
\label{trunc-cone}
R(b, \varphi ,r, s):=\big(B^n(b, r)\setminus \overline{B}^n(b, s)\big)\cap C(b,\varphi)
\end{equation}
for the open truncated cone at $b$.

\subsection*{Quasiregular mappings.}

A continuous function $f\colon G\to\mathbb{R}^n$, $n\geq2$, of a domain $G$ in
$\mathbb{R}^n$ is called {\it quasiregular} if $f$ is in the Sobolev
space $W^{1,n}_{\mathrm{loc}}(G)$ and there
 exists a constant $K$, $1\leq K<\infty$, such that the inequality
\begin{equation}
\label{dilatation}
|f'(x)|^n\leq KJ_f(x)
\end{equation}
holds a.e. in $G$. Here $f'(x)$ is the formal derivative and
$|f'(x)|=\max_{|h|=1}|f'(x)h|$.

The smallest $K\geq 1$ for which this inequality \eqref{dilatation} is true is called the outer dilatation of $f$ and
denoted by $K_O(f)$. If $f$ is quasiregular, then the smallest
$K\geq1$ for which the inequality
$$
J_f(x) \leq K l(f'(x))^n
$$
holds a.e. in $G$ is called the inner dilatation of $f$ and denoted by $K_I(f)$, where $l(f'(x))=\min_{|h|=1}|f'(x)h|$. The maximal
dilatation of $f$ is the number $K(f):=\max\{K_I(f),K_O(f)\}$. If $K(f)\leq K<\infty$, then $f$ is said to be $K$-quasiregular.

\subsection*{Hyperbolic type metrics}
Recall that the {\it hyperbolic distance} $\rho_{\Bn}$ in the unit ball $\Bn$, $n\ge 2$, also denoted by $\rho$, is defined by \cite[p. 55]{HKV}
\begin{equation}
\rho_{\Bn}(x,y) := \inf_\gamma \int_\gamma \frac{2|dz|}{1-|z|^2},
\end{equation}
where the infimum is taken over all rectifiable curves $\gamma$ joining $x$ and $y$ in $\Bn$.

The  hyperbolic distance in $\Bn$ can also be defined by using the following well-known formula:
\begin{equation}
\label{qr218}
\sinh^2\Big(\frac{1}{2}\rho(x,y)\Big)
= \frac{|x-y|^2}{(1-|x|^2)(1-|y|^2)},
\,\,\,x,y\in \Bn.
\end{equation}

A useful special case of this formula for $t\in(0,1)$ and $s\in(-t,t)$ is the formula \cite[2.17]{VuorinenBook}:
\begin{equation}
\label{hypest}
\rho(se_1,te_1) = \log\bigg(\frac{1+t}{1-t}\cdot \frac{1-s}{1+s}\bigg),
\end{equation}
where $e_1 = (1,0,\ldots,0)$ denotes the first unit basis vector of $\mathbb{R}^n$.

For a proper subdomain $G$ of $\IR^n$, define the {\it quasihyperbolic distance} as
$$k_G(x, y):=\inf\limits_{\gamma}\int_{\gamma}d(x)^{-1}\,ds,$$
where $\gamma$ runs through all rectifiable curves with $a, b\in \gamma$. It is clear that $k_G$ is a metric in $G$. For $x\in G$, $M>0$, we write $B_k(x, M):=\{z\in G: k_G(x, z)<M\}$.

%\subsection*{Distance ratio metrics}
For a domain $G$ in $\IR^n$, $G\neq \IR^n$, let $d(z)=d(z,\partial G)$ for $z\in G$. We define the {\it distance ratio metric} by
$$j_G(x, y):=\log\bigg(1+\frac{|x-y|}{\min\{d(x), d(y)\}}\bigg),$$
for $x, y\in G$. An elementary argument (see \cite[Lemma 4.6]{HKV}) shows that $j_G(x,y)$ is a metric on $G$.

If $A\subset G$ is non-empty, we define
$$ k_G(A) :=\sup\limits_{x, y\in A}k_G(x,y)\text{ and }
j_G(A):=\sup\limits_{x, y\in A}j_G(x,y).
$$
For more information about hyperbolic type metrics in $\mathbb{R}^n$, we refer to \cite{HKV}.

\subsection*{Angular sets and hyperbolic type metrics}

Next we will recall some elementary inequalities that will allow us to estimate the hyperbolic diameters of sets contained in conical domains.

\begin{prop}\label{j-estimate}
Let $\varphi \in (0, \pi/2)$, $b\in S^{n-1}$, $r\in (0, \cos (\varphi))$, and $a\in (0,1)$. %Set $R(b,\varphi, ar, r):=C(\varphi,b)\cap \big(\overline{B}^n(b, r)\setminus B^n(b, ar)\big)$.
Then
\begin{eqnarray*}
j_{\Bn}\big(R(b,\varphi, ar, r)\big)=\sup\limits_{x, y\in R(b,\varphi, ar, r)}\log\bigg(1+\frac{ |x-y|}{\min\{d(x), d(y)\}}\bigg) \leq \log\bigg(1+\frac{(2+ar)u(a, \varphi)}{a(2\cos(\varphi)-ar)}\bigg),
\end{eqnarray*}
where $u^2(a, \varphi)=(1+a)^2\tan^2(\varphi)+(1-a)^2$, where $R(b,\varphi, ar, r)$ is as in \eqref{trunc-cone}.
\end{prop}

\bpf
We first calculate an upper bound for $d(R(b,\varphi, ar, r))$.

Let $B$ be the smallest Euclidean ball containing $R(b,\varphi, ar, r)$, and denote by $O'$ and $r'$ the center and the radius of $B$, respectively. Applying the Law of Cosines yields
$$\left\{
\begin{aligned}
r'^2=&r^2+|O'b|^2-2r|O'b|\cos(\varphi),\\ \nonumber
r'^2=&a^2r^2+|O'b|^2-2ar|O'b|\cos(\varphi)
\end{aligned}
\right. $$
which imply that
\[
|O'b|=\frac{(1+a)r}{2\cos(\varphi)}.
\]
By applying the Pythagorean theorem, we obtain
\begin{eqnarray*}
r'^2&=&|O'b|^2\sin^2(\varphi)+\Big(\frac{r-ar}{2}\Big)^2\\ \nonumber
&=&\frac{r^2}{4}\Big((1+a)^2\tan^2(\varphi)+(1-a)^2\Big).
\end{eqnarray*}
We write $u(a, \varphi)^2=(1+a)^2\tan^2(\varphi)+(1-a)^2$, which gives  $r'=ru(a, \varphi)/2$.
Therefore, $d(R(b,\varphi, ar, r))\leq 2r'=ru(a, \varphi)$.

Now we calculate a lower bound of $\min\{d(x), d(y)\}$. Denote $\nu :=d(R(b,\varphi, ar, r), \partial \mathbb{\mathbb{B}}^n)$, and set $\nu=1-\xi$, where $$\xi^2=(1-ar\cos(\varphi))^2+(ar\sin\varphi)^2=1-2ar\cos(\varphi)+a^2r^2.$$ It follows that $1+\xi\leq 2+ar$ and $$\nu=1-\xi=\frac{1-\xi^2}{1+\xi}\geq \frac{ar(2\cos(\varphi)-ar)}{2+ar}.$$  \\ Therefore, $$j_{\Bn}\big(R(b,\varphi, ar, r)\big)\leq \log\bigg(1+\frac{(2+ar)u(a, \varphi)}{a(2\cos(\varphi)-ar)}\bigg).$$

\epf

\begin{rem}\label{HMAremark2}
From \cite[Lemma 4.9 and Remark 5.4]{HKV}, we know that 
$$
j_{\mathbb{B}^n}\leq k_{\mathbb{B}^n}\leq \rho_{\mathbb{B}^n}\leq 2j_{\mathbb{B}^n} \le 2\rho_{\mathbb{B}^n}.
$$
Therefore, we have
$$
k_{\mathbb{B}^n}\big(R(b,\varphi, ar, r)\big)=\sup\limits_{x, y\in R(b,\varphi, ar, r)} k_{\mathbb{B}^n}(x, y)\leq
 2\log\bigg(1+\frac{(2+ar)u(a, \varphi)}{a(2\cos(\varphi)-ar)}\bigg).
 $$
Because $r\in(0, \cos(\varphi))$ and $u(a, \varphi)^2\leq (1+a)^2(\tan^2(\varphi)+1)=(1+a)^2/\cos^2(\varphi)$, we also have 
$$
k_{\mathbb{B}^n}\big(R(b,\varphi, ar, r)\big)
\leq 2\log\bigg(1+\frac{(2+a\cos(\varphi))(1+a)}{a(2-a)\cos^2(\varphi)}\bigg).
$$
\end{rem}

\begin{lem}\label{HRMlemma1}
Let $b\in S^{n-1}$, $\varphi\in (0,\pi/2)$, $x \in C(b,\varphi)$, and $r:= |x-b|\in (0,\cos(\varphi))$. Then $$\rho\big(0,(1-r)b\big) \leq \rho(0,x)\leq s(r, \varphi),$$ where
$$s(r, \varphi):=\log\bigg(\frac{(2+r)^2}{r(2\cos(\varphi)-r)}\bigg).$$
 In particular, the line segment $L (r):=\big(0,(1-r)b\big)$ is contained in the hyperbolic ball $B_\rho(x,s(r,\varphi))$.
\end{lem}

\bpf
For $x\in \Bn$, we have $$\rho(0, x)=\log\bigg(\frac{1+|x|}{1-|x|}\bigg),$$
where the function
$$
t\mapsto \frac{1+t}{1-t}
$$
is increasing for $t\in (0,1)$. Furthermore, we have
$$\rho\big(0,(1-r)b\big) \leq \rho(0,x),$$
 because $$1-r\leq 1-r\cos(\varphi')\leq |x|,$$ where $\varphi'\leq \varphi$.

For any $x\in C(b, \phi)$,  let $\theta$ be the angle $\measuredangle (x,0,b)$, where $\theta\leq \varphi$. Then we know that $|x|^2=1+r^2-2r\cos(\theta)$. It follows that
 $|x|^2\leq 1+r^2-2r\cos(\varphi)$ and $|x|\leq 1+r$. Then we have
%\begin{eqnarray}
\begin{equation}
\label{HRMeq1}
\rho(0, x)=\log\bigg(\frac{1+|x|}{1-|x|}\bigg)
=\log\bigg(\frac{(1+|x|)^2}{1-|x|^2}\bigg)
\leq\log\bigg(\frac{(2+r)^2}{r(2\cos(\varphi)-r)}\bigg).
\end{equation}
%\end{eqnarray}
\epf

\subsection*{Modulus of a curve family} Let $\Gamma$ be a family of curves in $\IR^n$. The modulus of
$\Gamma$ is defined by \cite{vaisala}
$$M(\Gamma):=\inf\limits_{\rho}\int_{\IR^n} \rho^n dm,$$
where $m$ is the $n$-dimensional Lebesgue measure and the infimum is taken over all nonnegative Borel-functions $\rho: \IR^n\rightarrow  \IR\cup \{ \infty \}$ with $\int_{\gamma}\rho\, ds\geq 1$ for all $\gamma\in \Gamma$. If $E$, $F$, and $G$ are subsets of $\IR^n$, then we write $$\triangle (E, F; G):=\{\gamma: [0, c) \rightarrow G\},$$
where $\gamma$ is continuous and non-constant, $\gamma(0)\in E$, and $\gamma(t)\rightarrow F$ as  $t\rightarrow c$.

\subsection*{Condenser and its capacity} A pair $(A, C)$ is said to be a condenser if $A\subset \IR^n$ is open and $C\subset A$ is a compact non-empty set. The capacity of $(A, C)$ is defined by \cite[p. 150]{HKV}
$$\capacity (A, C):=\inf\limits_{u}\int_{\IR^n} |\nabla u|^ndm,$$
where $u$ runs through all $C_0^{\infty}(A)$ functions with $u(x)\geq 1$ for $x\in C$.

Equivalently, the capacity of the condenser $(A,C)$ can be defined by (see \cite[Theorem 9.6]{HKV})
$$\capacity (A, C):=M\big(\triangle(C, \partial A; A)\big)=M\big(\triangle (C, \partial A; \IR^n)\big).$$

For $E\subset \IR^n$, $x\in \IR^n$, and $r>0$, we write $$M(E, r, x):=M\big(\triangle(S^{n-1}(x, 2r), \overline{B}^n(x, r)\cap E; \IR^n)\big),$$
and define the lower and upper capacity densities of $E$ at $x$ by
$$\capacitydenL (E, x):=\liminf\limits_{r\rightarrow 0} M(E, r, x),$$
$$\capacitydenU (E, x):=\limsup\limits_{r\rightarrow 0} M(E, r, x).$$

\subsection*{Covering lemma}
Recall the following geometric lemma (see e.g. \cite[Lemma 3.2, p. 197]{landkof}):

\begin{lem}
\label{covering}
Let $n\ge 2$, and suppose that a set $F\subset \IR^n$ is covered by balls such that 
 each point $x\in F$ is a center of a ball $B_x$ of radius $r(x)>0$. If $\sup_{x\in F} r(x)< \infty$, then from
  the covering $\{ B_x : x\in F\}$ it is possible to choose a countable sub-covering $\{ B_{x_k}: k\in \mathbb{N}\}$ of
  $F$ so that the multiplicity of the covering is no larger than a number $c_0=c_0(n)$ depending only on the dimension $n$.
\end{lem}

\subsection*{Two-constants theorem and quasihyperbolic metric}

An important tool in this investigation is the following lemma of the third author, which is based on the two-constants theorem for quasiregular mappings \cite[Theorem 6.18, p. 190]{Ri2}.

\begin{lem}\cite[Lemma 2.22]{Vu2}\label{Lemma 2.22}
Let $f: G \rightarrow \Bn$ be a $K$-quasiregular mapping, $F\subset G$ compact, $\capacity(G, F)\geq \delta>0$, and $|f(x)|\leq \varepsilon_0$ for $x\in F$. Then there is a number $\beta\in (0,1)$ depending only on $n ,K, \delta$, and $d(F)/d(F, \partial G)$ such that
$$|f(y)|\leq \varepsilon_0^{\beta^{k(y, F)c_1+1}},$$
for $y\in G$, where $k(y, F):=\inf\{t>0: F\subset \overline{B_k(y, t)}\}$, $c_1:=(\log(1+\lambda_K))^{-1}$, and $\lambda_K\in (0,1/2)$ is a constant depending only on $n$ and $K$.
\end{lem}

\begin{rem}\label{constants}
From the proof of \cite[Lemma 2.22]{Vu2}, we know that
\begin{equation}
\label{beta-def}
\beta=\beta(\delta, n, K):=c_2\min\bigg\{ \Big(\frac{\delta}{N}\Big)^{1/(n-1)}, \Big(\frac{c_3}{2c_1}\Big)^{1/(n-1)}\bigg\},
\end{equation}
%$$\beta:=\min\big\{c_2(\delta/N)^{1/(n-1)}, c_2[(c_n/2)\log(1+\lambda_K)]^{1/(n-1)}\big\},$$
where $c_2=c_2(n,K)>0$ is a constant depending only on $n$ and $K$,
$$N=N\big(n,K,d(F)/d(F,\partial G)\big):=c_0\bigg(2+\frac{d(F)}{d(F, \partial G)\lambda_K}\bigg)^n,$$
$c_0=c_0(n)$ is as in Lemma \ref{covering}, and the constant $c_3=c_3(n)$  has an expression $c_3(n)=2^n b_n$, where $b_n$ is the constant from the spherical
cap inequality. The number $b_n>0$ is given by
$$
%\left\{\begin{array}{llllllll}
b_n=2^{1-2n}\omega_{n-2}\bigg(\int_0^{\pi/2}\sin^{\frac{2-n}{n-1}}t\,dt\bigg)^{1-n},\quad b_2=\frac{1}{2\pi},
%I_n=,
%\end{array}\right.
$$
cf. \cite[10.2 and 10.9.]{vaisala} and \cite[5.28. and
5.30.]{VuorinenBook}. 

In particular, we investigate the case where the quantity $\delta>0$ in \eqref{beta-def} is not fixed, and $\delta\to 0$ according to a given parameter.
For this case, observe that for sufficiently small values of $\delta$, i.e.,  $\delta < c_3 N/(2c_1)$, we have the lower bound
\begin{equation}
\label{beta-est}
\beta(\delta) \ge c_2\Big(\frac{\delta}{N}\Big)^{1/(n-1)}.
\end{equation}

If we apply Proposition \ref{j-estimate} and \cite[Lemma 2.17]{Vu2} with $\varphi\in(0,\pi/2)$, $F=R(b,\varphi, ar, r)$, $r\in (0,\cos(\varphi))$, and $G=\mathbb{B}^n$, then we have the upper bound
\[
N(r)=c_0\bigg(2+\frac{d(R(b,\varphi, ar, r))}{d(R(b,\varphi, ar, r), \partial \mathbb{B}^n)\lambda}\bigg)^n\\ \nonumber
\leq c_0\bigg(2+\frac{ru(a, \varphi)}{\nu \lambda}\bigg)^n\\ \nonumber
\leq c_0\gamma(a, \varphi, \lambda, n),
\]
where $u(a, \varphi)$ is as in Proposition \ref{j-estimate}, $\nu :=d(R(b,\varphi, ar, r), \partial \mathbb{\mathbb{B}}^n)$, and
$$
\gamma(a, \varphi, \lambda, n):=\bigg(2+\frac{(2+a\cos(\varphi))(1+a)}{a(2-a)\cos^2(\varphi)\lambda}\bigg)^n.
$$
We also have the upper bound
\[
N(r)\leq c_0\widetilde{\gamma}(a, r, \lambda, n),
\]
where
$$
\widetilde{\gamma}(a, r, \lambda, n):=\bigg(2+\frac{(2+ar)(1+a)}{a(2-a)r^2\lambda}\bigg)^n.
$$
\end{rem}

\subsection*{Boundary behavior}
Lemma \ref{Lemma 2.22} was used in \cite{Vu2} to prove the following result:

\begin{thm}\cite[Theorem 3.1]{Vu2}\label{Theorem 3.1}
Let $f: \Bn\rightarrow \Bn$ be a quasiregular mapping, $b\in S^{n-1}$, $\varphi \in (0, \pi/2)$, and let $E\subset C(b, \varphi)$. If $f(x)$ tends to $0$ when $x$ approaches $b$ through $E$, and if $\capacitydenL (E, b)=\delta >0$, then $f$ has an angular limit $0$ at $b$.
\end{thm}

The main goal of this investigation is to improve the above theorem by removing the assumption that the lower capacity density of $E$ is bounded from below by a constant $\delta>0$. We will replace this condition by a weaker assumption which permits the lower capacity density of $E$ to be $0$ at $b$ but, at the same time, requires that the function $f(x)$ tends to $0$ fast enough when $x\to b$ along the set $E$.

\section{Lindel\"of type results}

We are now ready to state the first of our main results. The next result is a refinement of Theorem \ref{Theorem 3.1}. Here the uniform lower bound $\delta>0$ is replaced
with the quantity $\delta(r)$, $r\in (0,1)$, which is allowed to approach zero as $r\to 0$.

\begin{thm}
\label{main1}
Let $f: \Bn\rightarrow \Bn$ be a K-quasiregular mapping, $b\in S^{n-1}$, $E\subset C(b, \varphi)$, $\varphi \in(0, \pi/2)$, and set $E_r:=E\cap R(b,\varphi, r, r/4)$. For $r\in (0,\cos(\varphi))$, let
$$
\delta(r):=\capacity\big(B^n(b, 2r), E_r\big)>0.
$$
 If $|f(x)|\leq \varepsilon(r)$ for $x\in E_r$, and 
 $$
 \bigg(\frac{\delta(r)}{\gamma_1(\varphi, n, K)}\bigg)^{\alpha_1(\varphi, n, K)}\log(\varepsilon(r))\rightarrow -\infty \text{ as }r \to 0,
 $$
where
$$\alpha_1(\varphi, n, K):=\bigg(\frac{24c_1}{\cos(\varphi)}+1\bigg)\Big(\frac{1}{n-1}\Big),$$
\[
\gamma_1(\varphi, n, K):=\bigg(\frac{14\lambda_K\cos^2(\varphi)+5\cos(\varphi)+40}{7\lambda_K\cos^2(\varphi)}\bigg)^n,
\]
and $c_1>0$ and $\lambda_K \in (0,1/2)$ are as in Lemma \ref{Lemma 2.22}, then $f$ has an angular limit $0$ at $b$.
\end{thm}

\bpf
%If $\delta(r)\ge \delta_0>0$ for all $r\in (0,1)$, the claim follows from Theorem \ref{Theorem 3.1}, and there is nothing to prove. Therefore, we may assume that $\delta(r)\to 0$ as $r \to 0$.
Divide the interval $(0,\cos(\varphi))$ into the sets
\[
I_0 := \big\{ r : \delta(r) \ge c_3 N(r)/(2c_1) \big\} \text{ and } I_1 := (0,\cos(\varphi))\setminus I_0.
\]
Let $r\in (0,\cos(\varphi))$. Applying Lemma \ref{Lemma 2.22}, we have
$$|f(x)|\leq \varepsilon(r)^{\beta(r)^{k_{\Bn}(x, \overline{E}_r)c_1+1}}$$
for $x\in R_r:= R(b,\varphi, r, r/4)$, where
% $$\beta:=\min\big\{C(\delta(r)/N)^{1/(n-1)}, C[(c_n/2)\log(1+\lambda_K)]^{1/(n-1)}\big\},$$
% where 
 $\beta(r) =\beta(\delta(r),n,K)$ is as in \eqref{beta-def}.
 
% $C>0$, $c_n$, and $N$ are as in Remark \ref{constants}.

By the proof of \cite[Theorem 3.1]{Vu2} with $\lambda=1/4$, we have $k_{\Bn}(R_r)\leq 24/\cos(\varphi)$. Therefore, $$|f(x)|\leq \varepsilon(r)^{\beta(r)^{\alpha(\varphi, n, K)}}$$ for $x\in R_r$, where
\[
\alpha(\varphi, n, K):=\frac{24c_1}{\cos(\varphi)}+1 >1.
\]

{\it Case I:}
If $0$ is an accumulation point of $I_0$ and $r\to 0$ in $I_0$, 
we note that the quantity $\beta(r)\equiv \beta_0>0$, where the quantity
$$\beta_0:=c_2\bigg(\frac{c_3}{2c_1}\bigg)^{{1}/(n-1)}.$$
Thus, it follows  immediately from the assumptions that $\varepsilon(r)\to 0$ and $|f(x)|\rightarrow 0$ when $x\rightarrow b$ in $\bigcup_{r\in I_0}R_r$.

{\it Case II:}
For the case $x\rightarrow b$ in $\bigcup_{r\in I_1}R_r$, we may use the estimate 
\eqref{beta-est} for $\beta(r)$. Because of the upper bound 
$N(r)\leq c_0\gamma_1(\varphi, n, K)$, we have 
$$
\beta(r)\geq c_2\bigg(\frac{\delta(r)}{N(r)}\bigg)^{{1}/(n-1)}\geq 
c_2\bigg(\frac{\delta(r)}{c_0\gamma_1(\varphi, n, K)}\bigg)^{{1}/(n-1)}.
$$
Therefore,
$$
\beta(r)^{\alpha(\varphi, n, K)}\geq c_2^{\alpha(\varphi, n, K)}\bigg(\frac{\delta(r)}{c_0\gamma_1(\varphi, n, K)}\bigg)^{\alpha_1(\varphi, n, K)}>0,
$$
 where
\[
\alpha_1(\varphi, n, K):=\frac{1}{n-1}\alpha(\varphi, n, K)
\]
and
\[
\gamma_1(\varphi, n, K):=\gamma (1/4, \varphi, \lambda_K, n)=\bigg(\frac{14\lambda_K\cos^2(\varphi)+5\cos(\varphi)+40}{7\lambda_K\cos^2(\varphi)}\bigg)^n.
\]
From the assumption that 
$$
\bigg(\frac{\delta(r)}{\gamma_1(\varphi, n, K)}\bigg)^{\alpha_1(\varphi, n, K)}\log(\varepsilon(r))\rightarrow -\infty,
$$
we have 
$$
\beta(r)^{\alpha (\varphi, n, K)}\log(\varepsilon(r))\rightarrow -\infty,
$$
which means that 
$$
|f(x)|\leq \varepsilon(r)^{\beta(r)^{\alpha(\varphi, n, K)}}\rightarrow 0.
$$

Cases I and II together show that $|f(x)|\rightarrow 0$ when $x\rightarrow b$ in $C(b, \varphi)$. Now, applying Theorem \ref{Theorem 3.1} with the set $E:=C(b, \varphi)$, which has positive lower capacity density at $b$, shows that $f$ has an angular limit $0$ at the boundary point $b$.
\epf

Besides angular limits, it is also possible to consider limits of some other kinds such as limits of the function along suitable types of tangential sets.

\begin{thm}
\label{main2}
Let $f: \Bn\rightarrow \Bn$ be a $K$-quasiregular mapping, $b\in S^{n-1}$,  let
 $\varphi:\ (0,1)\rightarrow (0, \pi/2)$ be a decreasing function with 
$\lim\limits_{r\rightarrow 0}\varphi(r)=\pi/2$ so that $r<\cos(\varphi(r))$, and let
$E\subset \Bn$ with
$$ \capacity(\Bn, E_r)\geq \delta>0\,, \quad E_r:=E\cap R\big(b,\varphi(r), r, r/2)\big),$$
 for  all $r \in (0,1)$.  If
 \begin{eqnarray}\label{equation of thm2}
 \gamma_3(\varphi(r), n, K)^{\alpha_2(\varphi(r), n, K)}\log\big(\varepsilon(r)\big)\rightarrow -\infty,  \text{ as }r \to 0,
 \end{eqnarray}
 where
% $$\alpha_2(\varphi(r), n, K)=2\log\bigg(1+\frac{4+\cos(\varphi(r))}{\cos^2 (\varphi(r))}\bigg)\bigg(\frac{1}{\log(1+\lambda_K)}\bigg)+1,$$
 $$\alpha_2(\varphi(r), n, K)=2c_1\log\bigg(1+\frac{4+\cos(\varphi(r))}{\cos^2 (\varphi(r))}\bigg)+1,$$
 $$\gamma_3(\varphi(r), n, K):=\bigg(\frac{2\lambda_K\cos^2(\varphi(r))+\cos(\varphi(r))+4}{\lambda_K\cos^2(\varphi(r))}\bigg)^{n/(1-n)},$$
 $c_1>0$ is as in Lemma \ref{Lemma 2.22}, and $|f(x)|\leq \varepsilon(r)$ for $x\in E_r$, then $f$ has a limit $0$ at b along the set $G_{b, \varphi}$ defined by
 $$G_{b, \varphi}:=\bigg\{x:\ \arccos\bigg(\frac{|xb-1|}{|x-b|}\bigg)<\varphi(|x-b|)\bigg\}.$$
\end{thm}
\begin{rem}
	Note that it follows from $(\ref{equation of thm2})$ that $\varepsilon(r)\rightarrow 0$ when $r\rightarrow 0$. However, the converse is not necessarily true, because we also have $\lim\limits_{r\rightarrow 0}\gamma_3(\varphi(r), n, K)^{\alpha_2(\varphi(r), n, K)}=0.$
\end{rem}
\bpf
 By Lemma \ref{Lemma 2.22}, we have $$|f(x)|\leq \varepsilon(r)^{\beta(r)^{k_{\Bn}(x, \overline{E}_r)c_1+1}}$$ on $R_{r}:=R(b, \varphi(r), r, r/2)$, where  $\beta(r) =\beta(r, \delta,n,K)$ is as in \eqref{beta-def}.
 %$$\beta:=\min\Big\{C\Big(\frac{\delta}{N}\Big)^{1/(n-1)}, C[(c_n/2)\log(1+\lambda_K)]^{1/(n-1)}\Big\}.$$

Applying Remarks \ref{HMAremark2} and \ref{constants} with $a=1/2$, we have 
 $$k_{\Bn}(R_r)\leq 2\log\bigg(1+\frac{4+\cos(\varphi(r))}{\cos^2 (\varphi(r))}\bigg),$$
and
$$N(r)\leq c_0\gamma(1/2, \varphi(r), \lambda_K, n):=c_0\gamma_2(\varphi(r), n, K),$$
where $$\gamma_2(\varphi(r), n, K):=\bigg(\frac{2\lambda_K\cos^2(\varphi(r))+\cos(\varphi(r))+4}{\lambda_K\cos^2(\varphi(r))}\bigg)^n.$$
Therefore, we have
$$|f(x)|\leq \varepsilon(r)^{\beta(r)^{\alpha_2(\varphi(r), n, K)}}$$
 for $x\in R_{r}$, where
$$\alpha_2(\varphi(r), n, K)=2c_1\log\bigg(1+\frac{4+\cos(\varphi(r))}{\cos^2 (\varphi(r))}\bigg)+1.$$

We will again consider two cases. Divide the interval $(0,1)$ into the sets
\[
I_0 := \big\{ r : \delta \ge c_3 N(r)/(2c_1) \big\} \text{ and } I_1 := (0,1)\setminus I_0,
\]

{\it Case I:} If $r\in I_0$, then $\beta\equiv\beta_0>0$, where 
$$\beta_0:=c_2\bigg(\frac{c_3}{2c_1}\bigg)^{{1}/(n-1)}.$$
In this case, we have 
$$
\lim\limits_{r\rightarrow 0}\alpha_2\big(\varphi(r), n, K\big)\log \big(\gamma_3(\varphi(r), n, K)\big)=-\infty,
$$
if $0$ is an accumulation point of the set of all such $r<\cos(\varphi(r))$ that $r\in I_0$. 

Therefore, 
$$
\lim\limits_{r\rightarrow 0}\gamma_3(\varphi(r), n, K)^{\alpha_2(\varphi(r), n, K)}=\lim\limits_{r\rightarrow 0}e^{\alpha_2(\varphi(r), n, K)\log (\gamma_3(\varphi(r), n, K))}=0.
$$
It follows from the assumption $(\ref{equation of thm2})$ that $\varepsilon(r)\rightarrow 0$ and $|f(x)|\rightarrow 0$ when $x\rightarrow b$ and $r=|x-b|\in I_0$.

{\it Case II:} 
If  $r\in I_1$, then 
\begin{eqnarray*}
\beta(r)&\geq& c_2\bigg(\frac{\delta}{c_0\gamma_2(\varphi(r), n, K)}\bigg)^{1/(n-1)}=c_2(\delta/c_0)^{1/(n-1)}\gamma_2(\varphi(r), n, K)^{1/(1-n)}\\ \nonumber
&:=&c_2(\delta/c_0)^{1/(n-1)}\gamma_3(\varphi(r), n, K).
\end{eqnarray*}
By our assumption $(\ref{equation of thm2})$,
we have 
$$\varepsilon(r)^{\beta(r)^{\alpha_2(\varphi(r), n, K)}}\rightarrow 0,$$
in $I_1$.

By combining these cases, it follows that $f$ has a limit $0$ at $b$ along the tangential set $G_{b, \varphi}$ defined by
$$G_{b, \varphi}:=\bigg\{x:\  \arccos\bigg(\frac{|xb-1|}{|x-b|}\bigg)<\varphi(|x-b|)\bigg\}.$$
\epf

\section{Capacity condition and Koebe arcs}
\label{koebe}

A classical result of P. Koebe states that if an analytic function 
$f: \B\rightarrow \B$ tends to zero along a sequence of arcs $(C_j)$ 
in the unit disk, which approach a sub-arc in the boundary,  such 
that for all $j$, $d(C_j)>\eta$ for  some $\eta>0$, then $f$ must 
be identically zero. Several refinements of this result are known in the literature. 
In particular,  D. C. Rung \cite{Rung} proved a Koebe type theorem 
for bounded analytic functions by weakening the assumption about the 
Euclidean diameters of the continua $(C_j)$. In his results the 
Euclidean diameters may tend to zero when $j\to \infty$, but the values of 
the function on the sequence  $(C_j)$ are required to converge toward 
a limit fast enough, depending on the hyperbolic diameters of $(C_j)$.

Rung's approach has been further studied by the second and the 
third authors in \cite{Rasila,Vu3}, where analogues of the 
aforementioned Koebe arcs were studied with assumptions made 
in terms of conformal moduli. Results of this type for planar 
harmonic mappings without quasiregularity assumption have also 
been investigated in \cite{bshouty1,samy3}.
Our next theorem further extends these results by using the 
technique from the previous section.

\begin{thm}
\label{main3}
Let $f: \Bn \to \Bn$ be a $K$-quasiregular mapping, $\varphi\in (0, \pi/2)$ and $b\in S^{n-1}$. Let $E\subset C(b, \varphi)$ and set $E_r:=E\cap R(b, \varphi, r, r/2)$ with $$\delta(r):=\capacity\big(B^n(b,2r), E_r\big)>0.$$
Assume that $|f(x)| \leq \varepsilon(r)$ for all $r\in (0,\cos(\varphi))$, $x\in E_r$, and let
$$
\bigg(\frac{\delta(r)}{\gamma_4(r, n, K)}\bigg)^{\alpha_3(r, \varphi, n, K)}\log(\varepsilon(r))\rightarrow -\infty \text{ as }r \to 0,
$$
where
$$\alpha_3(r,\varphi, n, K):=\bigg(c_1\log\bigg(\frac{(2+r)^2}{r(2\cos(\varphi)-r)}\bigg)+1\bigg)\Big(\frac{1}{n-1}\Big),
$$
$$\gamma_4(r, n, K):=\bigg(\frac{2\lambda_Kr^2+r+4}{\lambda_Kr^2}\bigg)^n,$$
and $c_1>0$ and $\lambda_K \in (0,1/2)$ are as in Lemma \ref{Lemma 2.22}. Then $f \equiv 0$.
\end{thm}

\bpf
Because $|f(x)| \leq \varepsilon(r)$ for $x \in E_r$, by Lemmas \ref{HRMlemma1} and \ref{Lemma 2.22}, we have
$$|f(y)| \leq \varepsilon(r)^{\beta(r)^{k_{\Bn}(y,\overline{E}_r)c_1+1}} \leq \varepsilon(r)^{\beta(r)^{s (r, \varphi)c_1+1}},$$
for $y\in B_{\rho}(x,s(r,\varphi))$, where $\beta(r) :=\beta(\delta(r),n,K)$ is as in \eqref{beta-def}, and
%$$\beta:=\min\{C(\delta(r)/N)^{1/(n-1)}, C[(c_n/2)\log(1+\lambda_K)]^{1/(n-1)}\},$$
$$s(r, \varphi):=\log\bigg(\frac{(2+r)^2}{r(2\cos(\varphi)-r)}\bigg).$$
In particular,
\begin{equation}
\label{f-est}
|f(y)| \leq \varepsilon(r)^{\beta(r) ^{s(r, \varphi)c_1+1}},
\end{equation}
for all $y \in L(1/2):=(0,b/2)$ and sufficiently small $r$. Note that $L(1/2)\subset B_{\rho}(x, s(r, \varphi))$ for small values of $r$, as $s(r, \varphi)\rightarrow \infty$ when $r\rightarrow 0$.

Divide the interval $(0,\cos(\varphi))$ into the sets
\[
I_0 := \big\{ r : \delta(r) \ge c_3 N(r)/(2c_1) \big\} \text{ and } I_1 := (0,\cos(\varphi))\setminus I_0,
\]
and study these cases separately.

{\it Case I:} For the case $r\in I_0$, if $0$ is an accumulation point of $I_0$ and $r\to 0$ in $I_0$, we note that $\beta(r)\equiv \beta_0>0$, where the quantity 
$$
\beta_0:=c_2\bigg(\frac{c_3}{2c_1}\bigg)^{{1}/(n-1)}.$$
Furthermore, for sufficiently small values of $r$,  the inequality 
$$\beta_0^{n-1}\geq \frac{\delta(r)}{\gamma_4(r, n, K)},$$
where
$$
\gamma_4(r, n, K):=\widetilde{\gamma}(1/2, r, \lambda_K, n)=\bigg(\frac{2\lambda_Kr^2+r+4}{\lambda_Kr^2}\bigg)^n,
$$
holds. Thus, we have
$$\beta_0^{s(r, \varphi)c_1+1}=\beta_0^{(n-1)\alpha_3(r, \varphi, n, K)}\geq \bigg(\frac{\delta(r)}{\gamma_4(r, n, K)}\bigg)^{\alpha_3(r, \varphi, n, K)}.$$
It follows  from the assumptions that for all $y\in L(1/2)$, the right hand side of  (\ref{f-est}) goes to $0$ as $r\rightarrow 0$ in $I_0$.

{\it Case II:} 
For the case $r\in I_1$, since $N(r)\leq c_0\gamma_4(r, n, K)$,  we have 
$$
\beta(r)\geq c_2\bigg(\frac{\delta(r)}{N(r)}\bigg)^{{1}/(n-1)}
\geq c_2\bigg(\frac{\delta(r)}{c_0\gamma_4(r, n, K)}\bigg)^{{1}/(n-1)}.
$$
Therefore, 
$$
\beta(r)^{s(r, \varphi)c_1+1}\geq c_2^{s(r, \varphi)c_1+1}\bigg(\frac{\delta(r)}{c_0\gamma_4(r, n, K)}\bigg)^{\alpha_3(r, \varphi, n, K)}>0,
$$
where $$\alpha_3(r, \varphi, n, K):=\frac{s(r, \varphi)c_1+1}{n-1}.$$
From the assumption that
$$
\bigg(\frac{\delta(r)}{\gamma_4(r, n, K)}\bigg)^{\alpha_3(r, \varphi, n, K)}\log\big(\varepsilon(r)\big)\rightarrow -\infty,
$$
we have
\begin{equation}
\label{eps-limit}
\beta(r)^{s(r, \varphi)c_1+1}\log(\varepsilon(r))\rightarrow -\infty.
\end{equation}
Combining \eqref{f-est} and \eqref{eps-limit}, we see that for all $y\in L(1/2)$, the right hand side of  (\ref{f-est}) approaches $0$ as $r\rightarrow 0$ in $I_1$.

Because the choice of $y\in L(1/2)$ is arbitrary, we have shown that the right hand side of  (\ref{f-est}) approaches $0$ as $r\rightarrow 0$ in $(0,\cos(\varphi))$. Therefore $f(y)\equiv 0$ for all $y \in L(1/2)$. It follows that $f$ cannot be a discrete function in $\Bn$, so it must be constant.
\epf

\subsection*{Acknowledgments}
The second author (Antti Rasila) was partly supported by NSF of China under the number 11971124 and NSF of Guangdong Province under the numbers 2021A1515010326 and 2024A1515010467, and Li Ka Shing Foundation under the number 2024LKSFG06. This collaboration was made possible by the visit of the third author to GTIIT in April 2023.

\subsection*{Declarations of interest} None.

\end{document}